
\baselineskip=14pt
\parskip=10pt
\def\halmos{\hbox{\vrule height0.15cm width0.01cm\vbox{\hrule height
  0.01cm width0.2cm \vskip0.15cm \hrule height 0.01cm width0.2cm}\vrule
  height0.15cm width 0.01cm}}

\magnification=\magstephalf

\def\1{{\overline{1}}}
\def\2{{\overline{2}}}
\parindent=0pt
\overfullrule=0in

\def\frac#1#2{{#1 \over #2}}

\bf
\centerline
{
Using the ``Freshman's Dream'' to  Prove  Combinatorial Congruences
}
\rm
\bigskip
\centerline
{\it By Moa APAGODU and Doron ZEILBERGER}

{\bf Abstract.}
In a recent beautiful but technical article, William Y.C. Chen, Qing-Hu Hou, 
and Doron Zeilberger  developed an algorithm for finding and proving congruence identities (modulo primes)
of indefinite sums of many combinatorial sequences, namely those (like the Catalan and Motzkin sequences)
that are expressible in terms of constant terms of powers of Laurent polynomials.
We first give a leisurely exposition of their elementary but brilliant approach, and then  extend it in two directions. 
The Laurent polynomials may be of several variables, 
and instead of single sums we have multiple sums. In fact we even combine these two generalizations!
We conclude with some {\bf super-challenges}.

[Added June 26, 2016: Roberto Tauraso [T], has pointed out that most of our conjectured super-congruences
are already known, except for super-congruence 6, for which he supplied a beautiful proof].

{\bf Introduction}

In a recent elegant article ([CHZ]) the following type of quantities were considered

$$
\left(\sum_{k=0}^{rp-1}a(k)\right)\,\,\hbox{mod}\,\,p \,\,,
$$

where 

$\bullet$  $a(k)$ is a combinatorial sequence, expressible as the constant term of a power of a Laurent polynomial of a {\bf single} variable
(for example, the central binomial coefficient ${2k \choose k}$ is the coefficient of $x^0$ in $(x+ \frac{1}{x})^{2k}$) .

$\bullet$ $r$ is a specific positive integer .

$\bullet$ $p$ is an arbitrary prime .

Let  $x \equiv_p y$ mean $x \, \equiv \, y \,\,(\hbox{mod} \,\, p)$, in other words, that $x-y$ is divisible by $p$.

The [CHZ] method, while ingenious, is very elementary! The main ``trick'' is:

{\bf The Freshman's Dream Identity} ([Wi]): $(a+b)^p \equiv_p a^p + b^p$ \quad .

Recall that the easy proof follows from the Binomial Theorem, and noting that ${p \choose k}$ is divisible by $p$ except when
$k=0$ and $k=p$. This also leads to one of the many proofs of the {\bf grandmother} of all congruences, {\bf Fermat's Little Theorem},
$a^p  \equiv_p a$,  by  starting with  $0^p  \equiv_p 0$,  and applying induction to $(a+1)^p \equiv_p a^p + 1^p$.

The second ingredient  in the [CHZ] method is even more elementary! It is:

{\bf Sum of a Geometric Series}:
$$
\sum_{i=0}^{n-1} z^i = \frac{z^n \, - \, 1}{z \, - \,1} \quad .
$$

The focus in the Chen-Hou-Zeilberger ([CHZ]) paper was both computer-algebra implementation, and proving a general theorem
about a wide class of sums.
Their paper is rather technical, so its beauty  is lost to a wider audience.
Hence the first purpose of the present article is to give a  leisurely introduction to their method, and
illustrate it with numerous illuminating examples. The second, main, purpose, however, is to extend the method
in {\bf two} directions. The summand $a(k)$, may be the constant term of a Laurent polynomial of {\bf several} variables,
and instead of a {\bf single} summation sign, we can have {\it multi-sums}. In fact we can combine these two!

{\bf Notation}

The constant term of a Laurent polynomial $P(x_1,x_2,\ldots,x_n)$, alias
the coefficient of $x_1^0x_2^0..x_n^0$, is denoted by $CT[P(x_1,x_2,\ldots,x_n)]$.
The general coefficient of 
$x_1^{m_1}x_2^{m_2}..x_n^{m_n}$ in $P(x_1,x_2,\ldots,x_n)$ is denoted by 
$COEFF_{[x_1^{m_1}x_2^{m_2}..x_n^{m_n}]}P(x_1,x_2,\ldots,x_n)$. For example, 
$$
CT\left[ \frac{1}{xy}+3+5xy-x^3+6y^2\right]=3 \quad , \quad
COEFF_{[xy]}\left[ \frac{1}{xy}+3+5xy+x^3+6y^2\right]=5.
$$

We use the symmetric representation of integers in $(-p/2, p/2]$ when reducing modulo a prime $p$. 
For example,  6 mod $5=1$ and 4 mod $5 =-1$. 

{\bf Review of the Chen-Hou-Zeilberger Single Variable Case}

In order to motivate our generalization, we will first review, in more detail than given in [CHZ], some
of their elegant results. 
Let's start with the {\bf Central Binomial Coefficients},  sequence  A000984 in the great OEIS ([Sl], https://oeis.org/A000984).

{\bf Proposition 1}. For any prime $p \geq 5$, we have
$$
\sum_{n=0}^{p-1} { 2n \choose n} \,\,
\equiv_p  \,\,
\cases{
1 ,& if \quad  $p \equiv 1 \, (mod \,\, 3)$ ;\cr
-1 ,& if \quad  $p \equiv 2 \, (mod \,\, 3)$ .\cr}   
$$

{\bf Proof}: Using the fact that 
$$
{ 2n \choose n} = CT\left[ \frac{(1+x)^{2n}}{x^n} \right] \quad ,
$$
and the Freshman's Dream identity, $(a+b)^p\equiv_p a^p+b^p$, we have
$$
\sum_{n=0}^{p-1} {2n \choose n} = \sum_{n=0}^{p-1} CT\left[ \left(\frac{(1+x)^{2n}}{x^n}\right)\right]
= \sum_{n=0}^{p-1} CT\left[ \left(2+x+\frac{1}{x}\right)^n\right]
$$
$$
= CT\left[ \frac{\left(2+x+\frac{1}{x}\right)^p-1}{2+x+\frac{1}{x}-1} \right]
\equiv_p 
CT \left[ \frac{2^p+x^{p}+\frac{1}{x^p}-1}{1+x+\frac{1}{x}}\right] \quad (\hbox{By Freshman's Dream})
$$
$$
\equiv_p CT \left[ \frac{2+x^{p}+\frac{1}{x^p}-1}{1+x+\frac{1}{x}}\right](\hbox{By Fermat's Little theorem})
$$
$$
= CT \left[ \frac{1+x^{p}+\frac{1}{x^p}}{1+x+\frac{1}{x}}\right]
\,= \,CT\left[ \frac{1+x^p+x^{2p}}{(1+x+x^2)x^{p-1}}\right]
= COEFF_{[x^{p-1}]} \left[ \frac{1}{1+x+x^2}\right]
$$
$$
= COEFF_{[x^{p-1}]} \left[ \frac{1-x}{1-x^3}\right]\,
= COEFF_{[x^{p}]} \left ( \, \sum_{i=0}^{\infty} x^{3i+1} \,\right ) \, + \, 
 COEFF_{[x^{p}]} \left ( \, \sum_{i=0}^{\infty} (-1) \cdot x^{3i+2} \, \right) \quad .
$$
The result follows from extracting the coefficient of $x^p$ in the above geometrical series.

{\bf Proposition 1'}.
$$
\sum_{n=0}^{2p-1} { 2n \choose n} \,\,
\equiv_p \,\,
\cases{
3 ,& if \quad  $p \equiv 1 \, (mod \,\, 3)$ ;\cr
-3 ,& if \quad  $p \equiv 2 \, (mod \,\, 3)$ .\cr}  
$$

{\bf Proof}: 
$$
\sum_{n=0}^{2p-1} {2n \choose n} = \sum_{n=0}^{2p-1} CT \left[ \left(  2+x+\frac{1}{x} \right)^n\right]
=CT \left[ \frac{\left(2+x+\frac{1}{x}\right)^{2p}-1}{2+x+\frac{1}{x}-1}\right]
$$
$$
=CT\left[ \frac{\left(6+4x+\frac{4}{x}+x^2+\frac{1}{x^2}\right)^{p}-1}{2+x+\frac{1}{x}-1}\right]
\,\equiv_p \, CT\left[ \frac{\left(6+4x^p+\frac{4}{x^p}+x^{2p}+\frac{1}{x^{2p}}\right)-1}{2+x+\frac{1}{x}-1}\right]
$$
$$
\,=\, COEFF_{[x^{2p-1}]} \left[  \frac{1+4x^p}{1+x+x^2}\right]
\,=\, COEFF_{[x^{2p-1}]} \left[ \frac{1}{1+x+x^2}\right]  +4 \cdot COEFF_{[x^{p-1}]} \left[ \frac{1}{1+x+x^2}\right] 
$$
$$
\, = \,COEFF_{[x^{2p-1}]} \left[ \frac{1-x}{1-x^3}  \right]  +4 \cdot COEFF_{[x^{p-1}]} \left[ \frac{1-x}{1-x^3}\right] 
$$
$$
= \, COEFF_{[x^{2p-1}]} \left[ \frac{1}{1-x^3}  \right]  + COEFF_{[x^{2p-1}]} \left[ \frac{-x}{1-x^3}  \right]  
+4 \cdot COEFF_{[x^{p-1}]} \left[ \frac{1}{1-x^3}\right]  +4 \cdot COEFF_{[x^{p-1}]} \left[ \frac{-x}{1-x^3}\right] \,\,
$$
$$
=COEFF_{[x^{2p}]} \left[ \sum_{i=0}^{\infty} x^{3i+1}  \right]  + COEFF_{[x^{2p}]} \left[  \sum_{i=0}^{\infty} (-1) \cdot x^{3i+2} \right]  
$$
$$
+4 \cdot COEFF_{[x^{p}]} \left[   \sum_{i=0}^{\infty} x^{3i+1} \right]  +4 \cdot COEFF_{[x^{p}]} \left[   \sum_{i=0}^{\infty} (-1) \cdot x^{3i+2} \right] \,\, .
$$
The result follows from extracting the coefficients of $x^{2p}$ in the first two geometrical series above,
and the coefficient of $x^p$ in the last two.

The same method (of [CHZ]) can be used to find the `mod $p$' of $\sum_{n=0}^{rp-1} { 2n \choose n}$ for {\bf any}
positive integer $r$. This leads to the following proposition,
whose somewhat tedious proof we omit.

{\bf Proposition 1''}.  For any prime $p \geq 5$, and any positive integer, $r$, 
$$
\sum_{n=0}^{rp-1} { 2n \choose n} \,\,
\equiv_p \,\,
\cases{
\alpha_r ,& if \quad  $p \equiv 1 \, (mod \,\, 3)$ ;\cr
-\alpha_r ,& if \quad  $p \equiv 2 \, (mod \,\, 3)$ ,\cr}   
$$
where
$$
\alpha_r \, = \, \sum_{n=0}^{r-1} {{2n} \choose {n}} \quad .
$$
For the record, here are the first ten terms of the integer sequence $\alpha_r$:
$$
[1,3,9,29,99,351,1275,4707,17577,66187] \quad .
$$
The sequence $\alpha_r$ is Sequence A6134 ([Sl],{\tt https://oeis.org/A006134}).
Note that $\alpha_r$ is the number of ways of tossing a coin $<2r$ times
and getting as many Heads as Tails.

The most {\it ubiquitous} sequence in combinatorics is sequence A000108 in the  OEIS \hfill\break
([Sl], {\tt https://oeis.org/A000108}, that according to Neil Sloane is the longest entry!), the
super-famous {\it Catalan Numbers}, $C_n:=\frac{(2n)!}{n!(n+1)!}$, that count
{\it zillions} of combinatorial families (see [St] for some of the more interesting ones).

{\bf Proposition 2}. Let $C_n$ be the Catalan Numbers,  then,  for every prime $p \geq 5$,
$$
\sum_{n=0}^{p-1}  C_n \,\,
\equiv_p \,\,
\cases{
1 ,& if \quad  $p \equiv 1 \, (mod \,\, 3)$ ;\cr
-2 ,& if \quad  $p \equiv 2 \, (mod \,\, 3)$ .\cr}   \quad .
$$

{\bf Proof}:  Since  $C_n={{2n} \choose {n}} - {{2n} \choose {n-1}}$, it
is readily seen that $C_n=CT[(1-x)(2+x+\frac{1}{x})^{n}]$.
We have
$$
\sum_{n=0}^{p-1}  C_n \,\,
= \sum_{n=0}^{p-1} CT \left[(1-x)\left(2+x+\frac{1}{x}\right)^n\right]
= CT \left[ \frac{(1-x)\left(\left(2+x+\frac{1}{x}\right)^p-1\right)}{2+x+\frac{1}{x}-1}\right]
$$
$$
\equiv_p CT\left[ \frac{(1-x)\left(\left(2+x^p+\frac{1}{x^p}\right)-1\right)}{2+x+\frac{1}{x}-1}\right]
\, = \, COEFF_{[x^{p-1}]} \left[  \frac{1-x}{1+x+x^2}\right] \, = \, COEFF_{[x^{p-1}]} \left[  \frac{(1-x)^2}{1-x^3}\right]
$$

$$
=  COEFF_{[x^{p}]} \left[  \frac{x}{1-x^3}\right]
+COEFF_{[x^{p}]} \left[  \frac{-2x^2}{1-x^3}\right]
+COEFF_{[x^{p}]} \left[  \frac{x^3}{1-x^3}\right]
$$
$$
=  COEFF_{[x^{p}]} \left[  \sum_{i=0}^{\infty} 1 \cdot x^{3i+1}  \right]
+COEFF_{[x^{p}]} \left[    \sum_{i=0}^{\infty} (-2) \cdot x^{3i+2}   \right]
+COEFF_{[x^{p}]} \left[   \sum_{i=0}^{\infty} 1 \cdot x^{3i+3} \right] \quad ,
$$
and the result follows from extracting the coefficient of $x^p$ from the first or second geometric series above
(note that we would never have to use the third geometrical series, since $p>3$).

The same method (of [CHZ]) can be used to find the `mod $p$' of $\sum_{n=0}^{rp-1} C_n$ for {\bf any}
{\it specific} positive integer $r$. In fact, one can keep $r$ general, but then the proof is
rather tedious,  and we will spare the readers (and ourselves, from typing it).

{\bf Proposition 2'}. Let $C_n$ be the Catalan Numbers,  then,  for any positive integer $r$, we have
$$
\sum_{n=0}^{rp-1}  C_n \,\,
\equiv_p \,\,
\cases{
\beta_r ,& if \quad  $p \equiv 1 \, (mod \,\, 3)$ ;\cr
-\gamma_r ,& if \quad  $p \equiv 2 \, (mod \,\, 3)$ ,\cr}   
$$
where

$$
\beta_r = \sum_{n=0}^{r-1} C_n \quad, \quad
\gamma_r = \sum_{n=0}^{r-1} (3n+2)C_n \quad .
$$

For the record, the first ten terms of the sequence of integer pairs $[\beta_r, -\gamma_r]$ are
$$
[ [1,-2], [2,-7], [4,-23], [9,-78], [23,-274], [65, -988], [197, -3628], [626,-13495], [2076,-50675], [6918, -191673]] \quad .
$$
We note that the sequence $\beta_r$ is sequence A014137 in the OEIS ([Sl], {\tt https://oeis.org/A014137})
but at this time of writing (June 9, 2016), the sequence $\gamma_r$ is not there (yet).

Not as famous as the Catalan numbers, but not exactly obscure, are the {\it Motzkin numbers}, $M_n$,
sequence A001006 in the OEIS ([Sl],  {\tt https://oeis.org/A001006}), that may be defined by the
constant term formula
$$
M_n= CT \left [ \displaystyle{(1-x^2)\left(1+x+\frac{1}{x}\right)^n} \right ] \quad .
$$

{\bf Proposition 3}.   Let $M_n$ be the Motzkin numbers, then for any prime $p \geq 3$, we have
$$
\sum_{n=0}^{p-1}  M_n  \,
\equiv_p  \,\,
\cases{
2 ,& if \quad  $p \equiv 1 \, (mod \,\, 4)$ ;\cr
-2 ,& if \quad  $p \equiv 3 \, (mod \,\, 4)$ .\cr}   
$$

{\bf Proof}: 
$$
\sum_{n=0}^{p-1} M_n = \sum_{n=0}^{p-1} CT \left[(1-x^2)\left(1+x+\frac{1}{x}\right)^n\right]
= CT \left[ \frac{(1-x^2)\left(\left(1+x+\frac{1}{x}\right)^p-1\right)}{1+x+\frac{1}{x}-1}\right]
$$
$$
\, \equiv_p  \, CT \left[ \frac{(1-x^2)\left( 1+x^p+\frac{1}{x^p}-1 \right)}{1+x+\frac{1}{x}-1}\right]
\, = \,   CT \left[ \frac{(1-x^2)\left( x^p+\frac{1}{x^p} \right)}{x+\frac{1}{x}}\right]
\, = \,   CT \left[ \frac{x(1-x^2)\left( x^p+\frac{1}{x^p} \right)}{1+x^2}\right]
$$
$$
= COEFF_{[x^{p-1}]} \left[  \frac{1-x^2}{1+x^2}\right]
= COEFF_{[x^{p}]} \left[ \frac{x}{1+x^2}\right]  -  COEFF_{[x^{p}]} \left[ \frac{x^3}{1+x^2}\right] \quad .
$$
$$
= COEFF_{[x^{p}]} \left[ \sum_{i=0}^{\infty} (-1)^i x^{2i+1} \right]  
+ COEFF_{[x^{p}]} \left[   \sum_{i=0}^{\infty} (-1)^{i+1} x^{2i+3} \right] \quad ,
$$
and the result follows from extracting the coefficient of $x^p$ from the first and second geometric series above,
by noting that when $p \equiv \,1 \, (mod \,\, 4)$, $i$ is even in the first series, and odd in the second one, and
vice-versa when  $p \equiv \,3 \, (mod \,\, 4)$.

The same method, applied to a general $r$  yields

{\bf Proposition 3'}.  Let $M_n$ be the Motzkin numbers, and let
$p \geq 3$ be prime, then 
$$
\sum_{n=0}^{rp-1}  M_n  \,
\equiv_p 
\cases{
2\delta_r ,& if \quad  $p \equiv 1 \, (mod \,\, 4)$ ;\cr
-2\delta_r ,& if \quad  $p \equiv 3 \, (mod \,\, 4)$ ,\cr}   
$$
where $\delta_r$ is the sequence of partial sums of the central trinomial coefficients,
sequence A097893 in the OEIS([Sl], {\tt https://oeis.org/A097893})
whose generating function is
$$
\sum_{r=0}^{\infty} \delta_r x^r \, = \,
{\frac {1}{ \left( 1-x \right) \sqrt { \left( 1+x \right)  \left( 1-3\,x \right) }}} \quad .
$$

From the above  proofs, it is easy to observe that partial sums with upper summation limit of the form  $rp-1$, for $r >1$, 
can always be expressed in terms of the sum with upper summation limit $p-1$. 
This observation leads us to the following simplification of Theorem 2.1 in [CHZ].

{\bf Theorem 4}. Let $P(x)$  be a Laurent polynomial in $x$ and let $p$ be a prime.  Let $R(x)$ be the denominator, after clearing, of the expression 
$$
\frac{P(x^p)-1}{P(x)-1}. 
$$
Then, for any positive integer $r$ and Laurent polynomial $Q(x)$,
$$
\left( \sum_{n=0}^{rp-1} CT \left[P(x)^nQ(x)\right] \right) \hbox{mod} \,\, p\,\,,
$$

is congruent to a finite linear combination of shifts of the sequence of coefficients of the rational function $\displaystyle{\frac{1}{R(x)}}$.

{\bf Multi-Sums and Multi-Variables}

We now extend the Chen-Hou-Zeilberger method for discovery and proof of congruence theorems to multi-sums and multi-variables.

{\bf Proposition 5}. Let $p \geq 5$ be  prime, then
$$
\sum_{n=0}^{p-1}\sum_{m=0}^{p-1} { n+m \choose m}^2\,\,
\equiv_p \,\,
\cases{
1 ,& if \quad  $p \equiv 1 \, (mod \,\, 3)$ ;\cr
-1 ,& if \quad  $p \equiv 2 \, (mod \,\, 3)$ .\cr}  
$$

{\bf{Proof}}: Let  
$$
P(x,y)=\left(1+y\right)\left(1+\frac{1}{x}\right) \quad ,
$$
and 
$$
Q(x,y) = \left(1+x\right)\left(1+\frac{1}{y}\right) \quad ,
$$
then
$$
 { n+m \choose m}^2 =  { n+m \choose m} { n+m \choose n}=CT \left[ P(x,y)^n Q(x,y)^m\right]\,\,.
$$
We have
$$
\sum_{m=0}^{p-1}\sum_{n=0}^{p-1}  { m+n \choose m}^2  = 
\sum_{m=0}^{p-1}\sum_{n=0}^{p-1}  CT\left[ P(x,y) ^nQ(x,y)^m\right]
=\sum_{m=0}^{p-1} CT\left[\frac{(P(x,y)^{p}-1)Q(x,y)^m}{P(x,y)-1}\right]
$$
$$
=CT\left[\left(\frac{P(x,y)^{p}-1}{P(x,y)-1}\right) \left(\frac{Q(x,y)^{p}-1}{Q(x,y)-1}\right)\right] \quad .
$$

Using the Freshman's Dream, $(a+b)^p \, \equiv_p \, a^p+b^p$, we can pass to mod $p$ as above, and get
$$
\sum_{m=0}^{p-1}\sum_{n=0}^{p-1}  { m+n \choose m}^2  \equiv_pCT\left[\left(\frac{P(x^p,y^p)-1}{P(x,y)-1}\right ) 
\left(\frac{Q(x^p,y^p)-1}{Q(x,y)-1}\right)\right]
\equiv_p CT\left[\frac{(1+y^p+x^py^p)(1+x^p+x^py^p)}{(1+y+xy)(1+x+xy)x^{p-1}y^{p-1}}\right]
$$
$$
\equiv_p COEFF_{[x^{p-1}y^{p-1}]}\left[\frac{(1+y^p+x^py^p)(1+x^p+x^py^p)}{(1+y+xy)(1+x+xy)}\right]
\equiv_p COEFF_{[x^{p-1}y^{p-1}]}\left[\frac{1}{(1+y+xy)(1+x+xy)}\right] \quad .
$$
It is possible to show that the coefficient of $x^ny^n$ in the Maclaurin expansion of the rational function
$\frac{1}{(1+y+xy)(1+x+xy)}$ is $1$ when $n \equiv \,0 (mod \,\, 3)$,
$-1$ when $n \equiv \,1\,(mod \,\, 3)$, and $0$ when $n \equiv \,2\,(\,mod \,\, 3)$. 
One way is to do a partial-fraction decomposition, and extract the coefficient of $x^n$, getting
a certain expression in $y$ and $n$, and then extract the coefficient of $y^n$.
Another way is by using the Apagodu-Zeilberger algorithm ([AZ]), that outputs that 
the sequence of diagonal coefficients, let's call them $a(n)$, satisfy the recurrence equation 
$a(n+2)+a(n+1)+a(n)=0$,  with initial conditions $a(0)=1,a(1)=-1$. 

Based on computer calculations, we conjecture 

{\bf Conjecture 5'}. For any prime $p \geq 5$, and any pair of positive integers, $r$, $s$, we have
$$
\sum_{n=0}^{rp-1}\sum_{m=0}^{sp-1} { n+m \choose m}^2\,\,
\equiv_p \,\,
\cases{
\epsilon_{rs} ,& if \quad  $p \equiv 1 \, (mod \,\, 3)$ ;\cr
-\epsilon_{rs} ,& if \quad  $p \equiv 2 \, (mod \,\, 3)$ ,\cr}  
$$
where
$$
\epsilon_{rs} =\sum_{m=0}^{r-1} \sum_{n=0}^{s-1}  { n+m \choose m}^2 \quad .
$$

We finally consider partial sums of {\it trinomial coefficients}.

{\bf Proposition 6}.  Let $p>2 $ be prime, then we have
$$
\sum_{m_1=0}^{p-1}\sum_{m_2=0}^{p-1} \sum_{m_3=0}^{p-1} { m_1+m_2+m_3 \choose m_1,m_2,m_3}  \, \equiv_p 1 \, \,\, .
$$

{\bf Proof}: First observe that 
${ m_1+m_2+m_3 \choose m_1,m_2,m_3}= CT \left[\frac{(x+y+z)^{m_1+m_2+m_3}}{x^{m_1}y^{m_2}z^{m_3}}\right]$.

Hence
$$
\sum_{0 \leq m_1, m_2, m_3 \leq p-1}{ m_1+m_2+m_3 \choose m_1,m_2,m_3}  = \sum_{0 \leq m_1, m_2, m_3 \leq p-1} 
CT \left[(x+y+z)^{m_1+m_2+m_3}/(x^{m_1}y^{m_2}z^{m_3}) \right]
$$
$$
= CT\left[\sum_{0 \leq m_1, m_2, m_3 \leq p-1}\frac{ (x+y+z)^{m_1+m_2+m_3}}{x^{m_1}y^{m_2}z^{m_3}}\right]
$$
$$
= CT \left[
\left ( \sum_{m_1=0}^{p-1} \left(\frac{x+y+z}{x}\right)^{m_1} \right )
\left ( \sum_{m_2=0}^{p-1} \left(\frac{x+y+z}{x}\right)^{m_2} \right )
\left ( \sum_{m_3=0}^{p-1} \left(\frac{x+y+z}{x}\right)^{m_3} \right )
\right] 
$$
$$
=CT\left[  \frac{(\frac{x+y+z}{x})^p -1}{\frac{x+y+z}{x}-1} \, \cdot \,
 \frac{(\frac{x+y+z}{y})^p -1}{\frac{x+y+z}{y}-1} \, \cdot \,
\frac{(\frac{x+y+z}{z})^p -1}{\frac{x+y+z}{z}-1} \right]
$$
$$
=COEFF_{[x^{p-1}y^{p-1}z^{p-1}]} \left [ \frac{(x+y+z)^p -x^p}{y+z} \cdot \frac{(x+y+z)^p -y^p}{x+z}  \cdot
\frac{(x+y+z)^p -z^p}{x+y} \right ] \,.
$$
So far this is true for all $p$, not only $p$ prime. Now take it mod $p$ and get, 
using the Freshman's Dream in the form $ (x+y+z)^p \equiv_p x^p+y^p+z^p$, that 
$$
\sum_{m_1=0}^{p-1}\sum_{m_2=0}^{p-1} \sum_{m_3=0}^{p-1} { m_1+m_2+m_3 \choose m_1,m_2,m_3}  
\,\, \equiv_p \,\,
COEFF_{[x^{p-1}y^{p-1}z^{p-1}]} \left(\frac{y^p+z^p}{y+z} \cdot \frac{x^p+z^p}{x+z}\cdot \frac{x^p+y^p}{x+y}\right)
$$
$$
=COEFF_{[x^{p-1}y^{p-1}z^{p-1}]} 
\left ( \sum_{i=0}^{p-1} (-1)^i y^i z^{p-1-i} \right)
\left ( \sum_{j=0}^{p-1} (-1)^j z^j x^{p-1-j} \right)
\left ( \sum_{k=0}^{p-1} (-1)^k x^{k} y^{p-1-k} \right)
$$
$$
=COEFF_{[x^{p-1}y^{p-1}z^{p-1}]}  \left[
\sum_{0 \leq i,j,k<p} (-1)^{i+j+k} x^{p-1-j+k} y^{i+p-1-k} z^{p-1-i+j}
\right ] \quad .
$$

The only contributions to the coefficient of $x^{p-1}y^{p-1}z^{p-1}$  in the above triple sum come when
$i=j=k$, so the desired coefficient of $x^{p-1}y^{p-1}z^{p-1}$
is
$$
\sum_{i=0}^{p-1} (-1)^{3i}= \sum_{i=0}^{p-1} (-1)^{i}=(1-1+1-1+ \dots +1-1)+ 1 = 1 \quad . \quad \halmos
$$

With more effort, one can get the following generalization.

{\bf Proposition 6'}.  Let $p\geq 3$ be prime, and let $r,s,t$ be any positive integers, then
$$
\sum_{m_1=0}^{rp-1}\sum_{m_2=0}^{sp-1} \sum_{m_3=0}^{tp-1} { m_1+m_2+m_3 \choose m_1,m_2,m_3}  \, \equiv_p  \,\, \kappa_{rst} \quad,
$$
where
$$
\kappa_{rst} =
\sum_{m_1=0}^{r-1}\sum_{m_2=0}^{s-1} \sum_{m_3=0}^{t-1} { m_1+m_2+m_3 \choose m_1,m_2,m_3} \quad .
$$

The same method of proof used in Proposition 6 yields (with a little more effort) a {\it multinomial} generalization.

{\bf Proposition 7}.  Let $p \geq 3 $ be prime,  then
$$
\sum_{m_1=0}^{p-1} \dots \sum_{m_n=0}^{p-1} 
{ m_1+\dots m_n \choose m_1, \dots , m_n}  \, \equiv_p 1 \, \,\, .
$$

In fact, the following  also holds.

{\bf Proposition 7'}.  Let $p \geq 3$ be prime, and let $r_1, \dots, r_n$ be positive integers, then
$$
\sum_{m_1=0}^{r_1p-1} \dots \sum_{m_n=0}^{r_np-1} 
{ m_1+\dots m_n \choose m_1, \dots , m_n}  
\, \equiv_p  \kappa_{r_1 \dots r_n} \quad ,
$$
where
$$
 \kappa_{r_1 \dots r_n} \, =\,
 \sum_{m_1=0}^{r_1-1} \dots \sum_{m_n=0}^{r_n-1}  { m_1+\dots m_n \choose m_1, \dots , m_n}  
\, \,\, .
$$

{\bf Super Congruences}

If a congruence identity that is valid modulo a prime $p$, is also valid modulo $p^2$ (or better still, modulo $p^3$ and beyond)
then we have a {\it super-congruence}. The grandmother of all supercongruences is  {\bf Wolstenholme's Theorem}([Wo], see also [We]) that
asserts that
$$
{2p-1 \choose p-1} \, \equiv_{p^3}  1 \quad ,
$$
that improves on the weaker version ${2p-1 \choose p-1} \, \equiv_{p^2}  1$, first proved by Charles Babbage ([B]), better known
for more impressive innovations.

To our surprise, most (but not all!) the above congruences have super-congruence extensions.
The method of [CHZ], as it stands now, is not applicable, since the ``Freshman's Dream'' is {\it only} valid modulo $p$,
hence we have no clue how to prove them. We leave them as challenges to our readers.

{\bf Super-Conjecture 1}. For any prime $p \geq 5$ 
$$
\sum_{n=0}^{p-1} { 2n \choose n} \, \,
\equiv_{p^2}  \,\,
\cases{
1 ,& if \quad  $p \equiv 1 \, (mod \,\, 3)$ ;\cr
-1 ,& if \quad  $p \equiv 2 \, (mod \,\, 3)$ .\cr}   
$$
More generally

{\bf Super-Conjecture 1''}.  For any prime $p \geq 5$, and any positive integer, $r$, 
$$
\sum_{n=0}^{rp-1} { 2n \choose n} \,\,
\equiv_{p^2} \,\,
\cases{
\alpha_r ,& if \quad  $p \equiv 1 \, (mod \,\, 3)$ ;\cr
-\alpha_r ,& if \quad  $p \equiv 2 \, (mod \,\, 3)$ ,\cr}   
$$
where
$$
\alpha_r \, = \, \sum_{n=0}^{r-1} {{2n} \choose {n}} \quad .
$$

{\bf Update, June 16, 2016}: Dennis Stanton proved Super-Conjectures $1$ and $1'$. 
See his nice write-up  \hfill\break
{\tt http://www.math.rutgers.edu/\~{}zeilberg/mamarim/mamarimhtml/freshmanDennisStanton.pdf}.
His method should also work for Super-Conjectures $2$, and probably $2'$. 

{\bf Update, June 17, 2016}: Roberto Tauraso([T]) pointed out that Super-Conjecture $1$ is already known,
see [T] for references.

{\bf Super-Conjecture 2}. Let $C_n$ be the Catalan Numbers,  then,  for every prime $p \geq 5$,
$$
\sum_{n=0}^{p-1}  C_n \,\,
\equiv_{p^2} \,\,
\cases{
1 ,& if \quad  $p \equiv 1 \, (mod \,\, 3)$ ;\cr
-2 ,& if \quad  $p \equiv 2 \, (mod \,\, 3)$ .\cr}   \quad .
$$

{\bf Update, June 17, 2016}: Roberto Tauraso([T]) pointed out that Super-Conjecture $2$ is already known,
see [T] for references..

More generally,

{\bf Super-Conjecture 2'}. Let $C_n$ be the Catalan Numbers,  then,  for any positive integer $r$, we have
$$
\sum_{n=0}^{rp-1}  C_n \,\,
\equiv_{p^2} \,\,
\cases{
\beta_r ,& if \quad  $p \equiv 1 \, (mod \,\, 3)$ ;\cr
-\gamma_r ,& if \quad  $p \equiv 2 \, (mod \,\, 3)$ ,\cr}   
$$
where

$$
\beta_r = \sum_{n=0}^{r-1} C_n \quad, \quad
\gamma_r = \sum_{n=0}^{r-1} (3n+2)C_n \quad .
$$

We note that poor Motzkin does not seem to have a super-extension, but Proposition 5 sure does!

{\bf Super-Conjecture 5}. Let $p \geq 5$ be  prime, then
$$
\sum_{n=0}^{p-1}\sum_{m=0}^{p-1} { n+m \choose m}^2\,\,
\equiv_{p^2} \,\,
\cases{
1 ,& if \quad  $p \equiv 1 \, (mod \,\, 3)$ ;\cr
-1 ,& if \quad  $p \equiv 2 \, (mod \,\, 3)$ .\cr}  
$$

{\bf Update, June 17, 2016}: Roberto Tauraso([T]) pointed out that Super-Conjecture $5$ follows easily
from  what we called Super-Conjecture $1$, see [T] for details.

More generally

{\bf Super-Conjecture 5'}. For any prime $p \geq 5$, and any pair of positive integers, $r$, $s$, we have
$$
\sum_{n=0}^{rp-1}\sum_{m=0}^{sp-1} { n+m \choose m}^2\,\,
\equiv_{p^2} \,\,
\cases{
\epsilon_{rs} ,& if \quad  $p \equiv 1 \, (mod \,\, 3)$ ;\cr
-\epsilon_{rs} ,& if \quad  $p \equiv 2 \, (mod \,\, 3)$ ,\cr}  
$$
where
$$
\epsilon_{rs} =\sum_{m=0}^{r-1} \sum_{n=0}^{s-1}  { n+m \choose m}^2 \quad .
$$

The most pleasant surprise is that Propositions 6 and 6' can be ``upgraded'' to a cubic super-congruence,
i.e. it is still true modulo $p^3$.

{\bf Super-Conjecture 6}.  Let $p>2 $ be prime,  then we have
$$
\sum_{m_1=0}^{p-1}\sum_{m_2=0}^{p-1} \sum_{m_3=0}^{p-1} { m_1+m_2+m_3 \choose m_1,m_2,m_3}  \, \equiv_{p^3} 1 \, \,\, .
$$

{\bf Update, June 17, 2016}: Roberto Tauraso([T]) has supplied a beautiful  proof, highly recommended.

More generally

{\bf Super-Conjecture 6'}.  Let $p\geq 3$ be prime, and let $r,s,t$ be any positive integers, then
$$
\sum_{m_1=0}^{rp-1}\sum_{m_2=0}^{sp-1} \sum_{m_3=0}^{tp-1} { m_1+m_2+m_3 \choose m_1,m_2,m_3}  \, \equiv_{p^3}  \,\, \kappa_{rst} \quad,
$$
where
$$
\kappa_{rst} =
\sum_{m_1=0}^{r-1}\sum_{m_2=0}^{s-1} \sum_{m_3=0}^{t-1} { m_1+m_2+m_3 \choose m_1,m_2,m_3} \quad .
$$

To our bitter disappointment, Propositions 7 and 7', for $n \geq 4$ summation signs, do not have  super-upgrades.

{\bf Lots and Lots of Combinatorial Challenges}

Perhaps the nicest proof of Fermat's Little Theorem, $a^p \equiv_p a$, is Golomb's ([G]) combinatorial proof, 
that notes that $a^p$ is the number of (straight) necklaces with $p$ beads, using beads of $a$ different colors,
and hence $a^p-a$ is the number of such (straight) necklaces that are not all of the same color.
For any such necklace, all its $p$ circular rotations are distinct (since $p$ is prime), hence
the set of such necklaces can be divided into families, each of them with $p$ members, and hence
there are $(a^p-a)/p$ `circular' necklaces (without clasp), and this must be an integer.

Each and every quantity in the above propositions and conjectures counts a natural combinatorial family.
For example $\sum_{n=0}^{p-1} {2n \choose n}$ counts the number of binary sequences with the same number
of $0$'s and $1$'s whose length is less than $2p$. Can you find a member of this set that
when you remove it, and $p \, \equiv \, 1\, (mod \, 3)$, you can partition that set into
families each of them with exactly $p$ (or better still, for the super-congruence, $p^2$) members?, and when
$p \, \equiv \, 2\,(mod \, 3)$,  can you find two such members?

{\bf References}

[AZ] Moa Apagodu and Doron Zeilberger,
{\it Multi-Variable Zeilberger and Almkvist-Zeilberger Algorithms and the Sharpening of Wilf-Zeilberger Theory},
Adv. Appl. Math. {\bf 37} (2006), 139-152. \hfill\break
{\tt http://www.math.rutgers.edu/\~{}zeilberg/mamarim/mamarimhtml/multiZ.html} \quad .

[B] Charles Babbage, {\it Demonstration of a theorem relating to prime numbers}, Edinburgh
Philosophical J. {\bf 1}(1819), 46-49.

[CHZ]  William Y.C. Chen, Qing-Hu Hou, and Doron Zeilberger, 
{\it Automated Discovery and Proof of Congruence Theorems for Partial Sums of Combinatorial Sequences}, J. of Difference Equations and Applications
(on-line before print),  DOI:10.1080/10236198.2016.1142541. Volume and page tbd.
Available from \hfill\break
{\tt http://www.math.rutgers.edu/\~{}zeilberg/mamarim/mamarimhtml/ctcong.html} \quad .

[G] Solomon W. Golomb, {\it Combinatorial proof of Fermat's ``Little'' Theorem},
American Mathematical Monthly {\bf 63}(1956), 718.

[Sl] Neil J. A. Sloane, {\it The On-Line Encyclopedia of Integer Sequences}, {\tt http://oeis.org/} \quad .

[St] Richard P. Stanley, {\it ``Catalan Numbers''}, Cambridge University Press, 2015.

[T] Roberto Tauraso, {\it A (human)  proof of a triple binomial sum supercongruence}, 
{\tt http://arxiv.org/abs/1606.05543} \quad .

[We] Eric W. Weisstein, {\it Wolstenholme's Theorem}, From MathWorld--A Wolfram Web Resource.  \hfill\break
{\tt http://mathworld.wolfram.com/WolstenholmesTheorem.html} \quad .

[Wi] The Wikipedia Foundation, {\it The Freshman's Dream} \quad .

[Wo] Joseph Wolstenholme, {\it On certain properties of prime numbers}, Quart. J. Pure Appl. Math. {\bf 5} (1862), 35-39.

\bigskip
\hrule
\bigskip
Moa Apagodu, Department of Mathematics, Virginia Commonwealth University, Richmond, VA 23284, USA;
$\,$ mapagodu at vcu dot edu .
\bigskip
\hrule
\bigskip
Doron Zeilberger, Department of Mathematics, Rutgers University (New Brunswick), Hill Center-Busch Campus, 110 Frelinghuysen
Rd., Piscataway, NJ 08854-8019, USA; 
$\,$ zeilberg at math dot rutgers dot edu  .
\bigskip
\hrule

\bigskip
{\bf  First Written:  June 9, 2016} ; 
{\bf  This version:  June 26, 2016} .

\end